\documentclass[12pt,a4paper]{amsart}
\usepackage{amsfonts}
\usepackage{amssymb}
\usepackage{amsmath}
\usepackage{amsthm}
\usepackage{cite}
\usepackage{enumitem}
\usepackage{tikz}
\usepackage{subfigure}

\theoremstyle{plain}
\newtheorem{thm}{Theorem}

\newtheorem{cor}{Corollary}
\newtheorem{bem}{Remark}

\providecommand{\R}{\mathbb{R}}

\providecommand{\eps}{\varepsilon}

\renewcommand{\qed}{\hfill $\Box$}
\newcommand{\derivord}[2]{
\frac{d #1}{d #2}
}

\providecommand{\R}{\mathbb{R}}

\providecommand{\eps}{\varepsilon}

\providecommand{\dr}{\,dr}
\providecommand{\dt}{\,dt}

\begin{document}
\allowdisplaybreaks

 \parindent0mm

\allowdisplaybreaks

\title{Uniqueness results for semilinear elliptic systems on $\R^n$}

\author{Rainer Mandel}
\address{R. Mandel \hfill\break
Department of Mathematics, Karlsruhe Institute of Technology (KIT)\hfill\break
D-76128 Karlsruhe, Germany}
\email{Rainer.Mandel@kit.edu}
\date{16.05.2013}


\begin{abstract}
  In this paper we establish uniqueness criteria for positive radially symmetric finite energy solutions of
  semilinear elliptic systems of the form
  \begin{align*}
    \begin{aligned}
    - \Delta u &= f(|x|,u,v)\quad\text{in }\R^n,\\
    - \Delta v &= f(|x|,v,u)\quad\text{in }\R^n.
    \end{aligned}  
  \end{align*}
  As an application we consider the following nonlinear Schr\"odinger system  
  \begin{align*}
    \begin{aligned}
    - \Delta u + u &= u^{2q-1} + b u^{q-1}v^q\quad\text{in }\R^n,\\
    - \Delta v + v &= v^{2q-1} + b v^{q-1}u^q \quad\text{in }\R^n.
    \end{aligned}  
  \end{align*}
  for $b>0$ and exponents $q$ which satisfy $1<q<\infty$ in case $n\in\{1,2\}$ and
  $1<q<\frac{n}{n-2}$ in case $n\geq 3$. Generalizing the results of Wei and Yao dealing with the case $q=2$
  we find new sufficient conditions and necessary conditions on $b,q,n$ such that precisely one
  positive solution exists. Our results dealing with the special case $n=1$ are optimal.
\end{abstract}

\maketitle

  \section{Introduction}
  
  In the paper \cite{WeiYao_Uniqueness_of_positive} Wei and Yao investigated whether positive finite
  energy solutions of the following nonlinear Schr\"odinger system
  \begin{align}\label{Eind Gl DGL Wei Yao}
    \begin{aligned}
     -\Delta u + u &= u^3 + buv^2 \quad\text{in }\R^n,\\
    - \Delta v + v &= v^3 + bvu^2 \quad\text{in }\R^n
    \end{aligned}
  \end{align}
  are uniquely determined when $n\in\{1,2,3\}$. One may summarize their main results as follows: If $b>1$ or
  $n=1,0<b<1$ then \eqref{Eind Gl DGL Wei Yao} has a unique positive solution $(u,v)$ and this
  solution satisfies $u=v$. The aim of our paper is to generalize the methods used in
  \cite{WeiYao_Uniqueness_of_positive} in order to prove uniqueness results for positive finite energy
  solutions of more general semilinear elliptic systems of the type
  \begin{align} \label{Gl PDGL f}
    \begin{aligned}
    - \Delta u &= f(|x|,u,v)\quad\text{in }\R^n,\\
    - \Delta v &= f(|x|,v,u)\quad\text{in }\R^n.
    \end{aligned}  
  \end{align}
  In some cases it is known that all such solutions  are radially
  symmetric. For instance, if $f(r,u,v)=f(u,v)$ is continuously differentiable on $\R_{\geq 0}^2$ and
  $f(0,0)=0, (f_u+f_v)(0,0)<0$ and $f_v\geq 0$ in every point of $\R_{\geq 0}^2$ then the symmetry results of Busca, Sirakov
  \cite{BusSir_Symmetry_results_for} and Ikoma \cite{Iko_Uniqueness_of_positive} imply that every
  positive finite energy solution of \eqref{Gl PDGL f} is radially symmetric and radially decreasing.
  Especially in these cases it is justified to consider the ODE version of \eqref{Gl PDGL f} given by
  \begin{align} \label{Eind Gl Thm 1 DGL}
    \begin{gathered}
      - u'' -\frac{n-1}{r}u' =  f(r,u,v)  \qquad\text{in }\R_{>0},  \\
      - v'' -\frac{n-1}{r}v'  = f(r,v,u)  \qquad\text{in }\R_{>0}, \\
      u'(0)=v'(0)=0, \quad  u(r),v(r)\to 0 \quad (r\to\infty).
    \end{gathered}
  \end{align}  
  In Theorem 1 and Theorem 2 we show that under suitable assumptions on $f$ every positive solution
  $(u,v)$ of \eqref{Eind Gl Thm 1 DGL} satisfies $u=v$ so that the uniqueness of the positive
  solution of \eqref{Eind Gl Thm 1 DGL} may be deduced from results for the corresponding scalar
  equation (cf.
  \cite{Jan_Uniqueness_of_positive},\cite{Kwo_Uniqueness_of_positive},\cite{PelSer_Uniqueness_of_positive}).
  In Corollary~\ref{Eind Cor} we apply these results to the family of nonlinearities given by
  \begin{equation} \label{Gl spezielles f}
    f(r,z_1,z_2)=-z_1+z_1^{2q-1}+bz_1^{q-1}z_2^q \qquad (r,z_1,z_2>0)
  \end{equation}  
  where $b$ is positive and $q$ satisfies
  \begin{equation} \label{Gl superlinear subcritical}
    1<q<\infty\quad \text{if }n\in\{1,2\}\qquad\text{and}\qquad 1<q<\frac{n}{n-2} \quad\text{if }n\geq 3.
  \end{equation}
  In the special case $q=2$ we reproduce the uniqueness criteria obtained in
  \cite{WeiYao_Uniqueness_of_positive}. Our results dealing with the cases $1<q<2$ and $q>2$ seem to be new.   
  
  \medskip
  
  In our first result dealing with \eqref{Eind Gl Thm 1 DGL} we need the following assumptions on $f$:
  \begin{itemize}
    \item[$(A1)$] The function $f:\R_{\geq 0}^3\to \R$ is continuous, locally Lipschitz continuous on
    $\R_{>0}^3$ and there are $p>1$ and $C>0$ such that
    $$
      |f(r,z_1,z_2)|\leq C (|z_1|+|z_2|+|z_1|^p+|z_2|^p)
      \quad\text{for all }r,z_1,z_2 \geq 0.
    $$
    \item[$(A2)$] There are positive numbers $m,R,\eps>0$ such that the following two inequalities
    hold for all $z_1,z_2\in (0,\eps),z_1\neq z_2$ and all $r>R$:
    \begin{align*}
      \frac{f(r,z_1,z_2)-f(r,z_2,z_1)}{z_1-z_2} \leq -m,\qquad
      \frac{f(r,z_1,z_2)+f(r,z_2,z_1)}{z_1+z_2} \leq -m.
    \end{align*}
    \item[$(A3)$] $f(r,z_1,z_2)z_2<f(r,z_2,z_1)z_1$ for all $z_1>z_2>0$ and all $r>0$.
  \end{itemize}  

  \begin{thm} \label{Eind Thm 1}
    Let $(u,v)$ be a positive solution of \eqref{Eind Gl Thm 1 DGL} where $f$ satisfies 
    $(A1),(A2),(A3)$. Then we have $u=v$ and there is $c>0$ such that $u(r)+|u'(r)|\leq
    c\, e^{-\sqrt{m}r}$ for all $r>0$ where $m$ is the positive number from $(A2)$.
  \end{thm}
 
  For nonlinearities that satisfy the opposite inequality in (A3) we find a uniqueness result
  only under additional assumptions. First we have to assume that the solution $(u,v)$ is
  decreasing, i.e., both $u$ and $v$ are decreasing. 
  More importantly we have to require that $f$ can be written in the form $f(r,u,v)= g(r,u) +
  h(r,u,v)v$ where $h(r,u,v)=h(r,v,u)>0$ and $g$ satisfies the differential inequality
  $rg_r(r,u)+(2n-2)g(r,u)\leq 0$. Unfortunately, this condition turns out to be quite restrictive for space
  dimensions $n\geq 2$ so that it would be desirable to find out how this assumption can be avoided.
  
  \begin{itemize} 
	\item[$(A3')$] We have $f(r,z_1,z_2)z_2> f(r,z_2,z_1)z_1$ for all $z_1>z_2>0$ and all $r>0$.
	\item[$(A4)$] There are functions $g\in C^1(\R_{>0}^2,\R)$ and $h\in C(\R_{>0}^3,\R)$
	such that for all $z_1,z_2,r>0$ we have
	$$
	  f(r,z_1,z_2)= g(r,z_1) + h(r,z_1,z_2)z_2
	$$
	where $h(r,z_1,z_2)=h(r,z_2,z_1)>0$ and $rg_r(r,z_1)+(2n-2)g(r,z_1)\leq 0$.
  \end{itemize}
  
  \begin{thm} \label{Eind Thm 2}
    Let $(u,v)$ be a positive decreasing solution of \eqref{Eind Gl Thm 1 DGL} where $f$
    satisfies the assumptions $(A1),(A2),(A3'),(A4)$. Then we have $u=v$ and there is $c>0$ such
    that ${u(r)+|u'(r)|\leq c\, e^{-\sqrt{m}r}}$ for all $r>0$ where $m$ is the positive number from
    $(A2)$.
  \end{thm}
  
  \begin{bem} \label{Bem Symmetrie}   
    Symmetry results based on the moving plane method imply the monotonicity of both component functions, see
    e.g. \cite{BusSir_Symmetry_results_for},\cite{Iko_Uniqueness_of_positive}. Hence, if for instance the
    nonlinearity $f=f(u,v)\in C^1(\R_{\geq 0}^2)$ satisfies the assumptions $(A1),(A2),(A3'),(A4)$ with
    differentiable functions $g,h$ such that $g_u(0)<0$ then a short calculation shows    
    $f(0,0)=0, (f_u+f_v)(0,0)<0$ and $f_v\geq 0$. In case $n\geq 1$  Theorem 1 in
    \cite{BusSir_Symmetry_results_for} implies that every positive solution of \eqref{Gl PDGL f} is
    radially symmetric and radially decreasing so that the monotonicity assumption of Theorem~\ref{Eind Thm
    2} is automatically satisfied. In case $n=1$ this follows from Theorem 4.1 in 
    \cite{Iko_Uniqueness_of_positive}.
  \end{bem}
  
  In Corollary \ref{Eind Cor} we apply the above theorems to the nonlinear Schr\"odinger system
  \begin{align}\label{Eind Gl DGL}
     \begin{gathered}
      - u'' -\frac{n-1}{r}u' + u = u^{2q-1} + bu^{q-1}v^q \qquad\text{in }\R_{>0}, \\
      - v'' -\frac{n-1}{r}v' + v = v^{2q-1} + bv^{q-1}u^q  \qquad\text{in }\R_{>0}, \\
      u'(0)=v'(0)=0, \quad u,v>0,\quad u(r),v(r)\to 0 \quad  (r\to\infty).
    \end{gathered}
  \end{align}
  Several authors proved the existence of positive solutions of \eqref{Eind Gl DGL} under suitable
  assumptions on $b$ and $q$, see for example \cite{MaMoPe_Positive_solutions_for} or
  \cite{Sir_Least_energy_solutions}. As mentioned before Wei and Yao proved in
  \cite{WeiYao_Uniqueness_of_positive} that in case $n\in\{1,2,3\},q=2,b>1$ or $n=1,q=2,0<b<1$ there is
  precisely one positive solution of \eqref{Eind Gl DGL} which is given by $(1+b)^{-\frac{1}{2}}(w_0,w_0)$
  where $w_0$ is the unique positive function satisfying $-\Delta w_0+w_0=w_0^3$ in $\R^n$. In Corollary \ref{Eind
  Cor} we generalize this result to all exponents $q$ satisfying \eqref{Gl superlinear subcritical}. For
  notational convenience we set
  \begin{equation} \label{Gl Def ubvb}
    (u_b,v_b) := (1+b)^{-\frac{1}{2q-2}}(u_0,u_0)
  \end{equation}
  where $u_0\in H^1(\R^n)$ is the unique positive function satisfying $-\Delta u_0+u_0=u_0^{2q-1}$ in
  $\R^n$.

  \begin{cor} \label{Eind Cor}
    Let the parameters $n,b,q$ satisfy \eqref{Gl superlinear subcritical} and
    \begin{itemize} 
      \item[(i)] $1<q<2,b\geq q-1$\quad or\quad $q=2,b>1$ \qquad or
      \item[(ii)] $n=1,q=2,0< b<1$\quad or \quad $n=1,q>2,0< b\leq q-1$.
    \end{itemize}  
    Then the function $(u_b,v_b)$ given by \eqref{Gl Def ubvb} is the unique positive solution of
    \eqref{Eind Gl DGL}.
  \end{cor}
  
  \medskip
  
  \begin{bem} ~
  \begin{itemize} 
    \item[(i)] In case $n=1$ the uniqueness results from Corollary \ref{Eind Cor} are optimal for positive
    $b$. Indeed, if $1<q<2,0<b<q-1$ or $q>2,b>q-1$ and if $k_b$ denotes the unique zero of
    $1+bk^q-bk^{q-2}-k^{2q-2}=0$ in $(0,1)$ then the functions 
    $$
      (u_0,u_0),\; (\mu_bu_0,\mu_bk_b u_0),\;(\mu_bk_b u_0,\mu_bu_0)\quad
      \text{where } \mu_b=(1+bk_b^q)^{-\frac{1}{2q-2}} 
    $$
    define three different positive solutions of \eqref{Eind Gl DGL}. Moreover, in case $q=2, b=1$ there is
    a family of solutions given by $(\cos(\theta) u_0,\sin(\theta) u_0)$ where $\theta\in
    (0,\frac{\pi}{2})$ and it is even known that there are no other positive solutions in this case, cf.
    Theorem 1.2 in \cite{WeiYao_Uniqueness_of_positive}. 
    \item[(ii)] As in \cite{WeiYao_Uniqueness_of_positive} we have to leave open whether part (ii)
    is true for space dimensions $n\geq 2$. 
    \item[(iii)] The uniqueness issue for \eqref{Eind Gl DGL} is much more difficult when negative
    coupling parameters $b$ are considered. At least for $q=2$ and $b$ sufficiently close to $-1$
    the bifurcation results of Bartsch, Dancer and Wang \cite{BaDaWa_A_Liouville_Theorem} show that a large
    variety of positive solutions exists so that uniqueness cannot hold in this case.
    \item[(iv)] The uniqueness result found by Ma, Zhao (cf. Theorem 2 in
    \cite{MaZh_Uniqueness_of_ground}) can be derived from our Theorem \ref{Eind Thm 1}. This
    statement will be proved in the Appendix.
  \end{itemize}
  \end{bem}  
    
  Let us finally compare our results with those of Quittner and Souplet \cite{QuiSou_Symmetry_of_components}. 
  For certain nonlinearities $f$ independent of $r$ with
  \begin{equation} \label{Gl Bed QuiSou}
    (z_1-z_2)(f(z_1,z_2)-f(z_2,z_1))\leq 0 \quad\text{for all }z_1\geq z_2>0
  \end{equation} 
  they prove that every positive solution $(u,v)$ of \eqref{Gl PDGL f} satisfies $u=v$. For instance, their
  results apply to the particular nonlinearities $f(u,v)=u^rv^p$ with $0\leq r,p\leq \frac{n}{n-2},n\geq 3$
  or $f(u,v)=-\lambda u^3+uv^2$ in case $0<\lambda<\frac{1}{n-2},n\geq 4$, see Theorem 1.1 and Theorem 1.6 
  in~\cite{QuiSou_Symmetry_of_components}. Since these nonlinearities do not
  meet the assumption $(A2)$ our theorems can not be used to reproduce their results. 
  Vice versa, our main example \eqref{Gl spezielles f} does not satisfy the inequality \eqref{Gl
  Bed QuiSou} so that Corollary 1 cannot be derived from the results in \cite{QuiSou_Symmetry_of_components}.
  As a consequence, our results may be regarded as being independent from the ones proved in
  \cite{QuiSou_Symmetry_of_components}. There are, however, semilinear systems of the type
  \eqref{Gl PDGL f} where $u=v$ can be derived both from the results contained in
  \cite{QuiSou_Symmetry_of_components} and from Theorem \ref{Eind Thm 1}.
  One example for such a system is given by the nonlinearity $f(u,v)=-u-u^{2q-1}+\beta u^{q-1}v^q$ with
  $\beta>0$ which has been considered in \cite{MaZh_Uniqueness_of_ground}.

  \section{Proof of Theorem \ref{Eind Thm 1}}
    
    For a positive solution $(u,v)$ of \eqref{Eind Gl Thm 1 DGL} we set 
    $$
      w_1:= u - v,\qquad w_2:= u+v.
    $$ 
    We assume $w_1\neq 0$ and our aim is to lead this assumption to a contradiction. The functions $w_1,w_2$
    satisfy the differential equations
     \begin{equation} \label{Eind Thm 1 DGL w1w2}
      -w_1''-\frac{n-1}{r}w_1 + c_1(r) w_1 = 0,\qquad
      -w_2''-\frac{n-1}{r}w_2 + c_2(r) w_2 = 0
    \end{equation}
    where the functions $c_1,c_2$ are given by
    \begin{align*}
      c_1(r) &:= -\frac{f(r,u(r),v(r))-f(r,v(r),u(r))}{u(r)-v(r)},\\
      c_2(r) &:= -\frac{f(r,u(r),v(r))+f(r,v(r),u(r))}{u(r)+v(r)}
    \end{align*}
    whenever $u(r)\neq v(r)$. We set $c_1(r):=c_2(r):=m$ in case $u(r)=v(r)$.
    From $u(r),v(r)\to 0$ as $r\to\infty$ and assumption (A2) we infer that there is $r_1>0$ sufficiently
    large such that 
    \begin{equation} \label{Eind Gl Absch c1 c2}
      c_1(r),c_2(r)\geq m \qquad\text{whenever }r\geq r_1.
    \end{equation}
    From \eqref{Eind Thm 1 DGL w1w2} and \eqref{Eind Gl Absch c1 c2} we conclude that the functions
    $w_1,w_2$ do not have any positive local maxima or negative local minima on $[r_1,\infty)$. Hence,
    both $w_1$ and $w_2$ are monotone on $[r_1,\infty)$ and we henceforth assume without loss of generality
    that $w_1$ is positive and decreasing on that interval. (If $w_1$ is negative and increasing one
    may interchange the roles of $u,v$.)
    
    \medskip
    
    Now we prove that $u(r),u'(r),v(r),v'(r)$ decay exponentially to zero as $r\to\infty$. This will ensure
    that certain integrals over unbounded domains are well-defined. From \eqref{Eind Thm 1 DGL w1w2}
    and $w_1'\leq 0$ on~$[r_1,\infty)$ we get
    $$
      (-{w_1'}^2+mw_1^2)'
      =2w_1'(-w_1''+m w_1) 
      \geq \frac{2(n-1)}{r} {w_1'}^2
      \geq 0 \quad\text{on }[r_1,\infty).
    $$
    As a consequence the function $-{w_1'}^2+mw_1^2$ is increasing on that interval. This implies
    $-{w_1'}^2+mw_1^2\leq 0$ and thus $w_1'+\sqrt{m}w_1\leq 0$ on $[r_1,\infty)$ so that the function
    $e^{-\sqrt{m}r}w_1$ is decreasing on $[r_1,\infty)$. A similar argument shows that the function
    $e^{-\sqrt{m}r}w_2$ is decreasing so that $0<u(r)+v(r)\leq c\,e^{-\sqrt{m}r}$ for all $r\geq r_1$ and some
    $c>0$. Integrating the differential equations 
    $$
      -(r^{n-1}u')'=r^{n-1}f(r,u,v),\qquad 
      -(r^{n-1}v')'=r^{n-1}f(r,v,u)
    $$
    and using $(A1)$ we finally obtain
    \begin{equation} \label{Eind Gl exp Abf}
      u(r)+v(r)+|u'(r)|+|v'(r)|\leq c\, e^{-\sqrt{m}r} \qquad (r>0)
	\end{equation}   
	for some $c>0$.
	
	\medskip
        
    From $u'(0)=v'(0)=0$, \eqref{Eind Gl exp Abf} and \eqref{Eind Gl Thm 1 DGL} we get
    \begin{align*}
	  0
      &= \int_0^\infty \Big( r^{n-1}(u'v-v'u) \Big)'\dr \\
      &= \int_0^\infty r^{n-1}\Big( \Big(u''+\frac{n-1}{r}u'\Big)v-\Big(v''+\frac{n-1}{r}v'\Big)u \Big)\dr \\
      &= \int_0^\infty r^{n-1} \big( -f(r,u,v)v+f(r,v,u)u  \big) \dr.
    \end{align*}
    These integrals are well-defined due to assumption (A1) and \eqref{Eind Gl exp Abf}. The
    assumption (A3) implies that $w_1=u-v$ has a zero in $(0,\infty)$. Since $w_1$ was shown to be
    positive and decreasing on $[r_1,\infty)$ the function $w_1$ also has a last zero $r_0$ which
    satisfies $r_0\leq r_1$. From $w_1>0$ on $(r_0,\infty)$ and the Picard-Lindel\"of Theorem we obtain
    $w_1'(r_0)>0$. We get
    \begin{align*}
      0
      &> - r_0^{n-1}u(r_0)w_1'(r_0) \\
      &= - r_0^{n-1} \big( u'(r_0)v(r_0)-v'(r_0)u(r_0) \big)  \\
      &= \int_{r_0}^\infty \Big(r^{n-1}(u'v-v'u) \Big)'\dr \\
      &= \int_{r_0}^\infty r^{n-1}\Big( \Big(u''+\frac{n-1}{r}u'\Big)v-\Big(v''+\frac{n-1}{r}v'\Big)u \Big)\dr
      \\
      &= \int_{r_0}^\infty r^{n-1} \big(-f(r,u,v)v+f(r,v,u)u  \big) \dr
    \end{align*}
    where the last integral is positive due to assumption $(A3)$ and $u>v$ on $(r_0,\infty)$. Therefore, we
    obtain a contradiction and the proof is finished.     \qed
    
    \begin{bem}
      The semilinear elliptic systems considered by Quittner and Souplet \cite{QuiSou_Symmetry_of_components}
      are of the more general form 
      \begin{align} \label{Gl System QuiSou}
        \begin{aligned}
        -\Delta u &= f_1(|x|,u,v) \quad\text{in }\R^n,\\
        -\Delta v &= f_2(|x|,v,u) \quad\text{in }\R^n.
        \end{aligned}
      \end{align}
      Textual modifications of the proof given above yield that the conclusion of Theorem \ref{Eind Thm 1}
      remains true provided both functions $f_1,f_2$ satisfy $(A1)$ as well as
      $f_1(r,z_1,z_2)z_2<f_2(r,z_2,z_1)z_1$ for all $z_1>z_2>0,r>0$ and provided there are positive
      numbers $m,R,\eps>0$ such that
      \begin{align*}
        \frac{f_1(r,z_1,z_2)-f_2(r,z_2,z_1)}{z_1-z_2} \leq -m,\qquad
        \frac{f_1(r,z_1,z_2)+f_2(r,z_2,z_1)}{z_1+z_2} \leq -m
      \end{align*}
      whenever $0<z_1,z_2<\eps,z_1\neq z_2$ and $r>R$.       
    \end{bem}

  \section{Proof of Theorem \ref{Eind Thm 2}}

    As in the proof of Theorem \ref{Eind Thm 1} we argue by contradiction. We assume that there
    is a positive decreasing solution $(u,v)$ of \eqref{Eind Gl Thm 1 DGL} that satisfies $u\neq v$. As above
    we obtain that the functions $u,v,u',v'$ satisfy the inequality \eqref{Eind Gl exp Abf} for
    some $c>0$ and that we may without loss of generality assume $(u-v)'(r_0)>0$ as well as $u-v>0$ on
    $(r_0,\infty)$ where $r_0$ denotes the last zero of $u-v$. We set $G(r,z):= \int_0^z g(r,t)\dt$.
    From the differential equation \eqref{Eind Gl Thm 1 DGL} and $(A4)$ we get
    \begin{align}\label{Gl Beweis Thm 2}
      \derivord{}{r} \Big( r^{2n-2}&\Big(-\frac{1}{2}u'^2+\frac{1}{2}v'^2-G(r,u)+G(r,v) \Big) \Big)
      \notag\\
       &= r^{2n-2}\Big( u' \Big(-u''-\frac{n-1}{r}u'-g(r,u) \Big) - v' \Big(
       -v''-\frac{n-1}{r}v'-g(r,v)\Big)\Big) \notag\\
       &\;  +r^{2n-2}(-G_r(r,u)+G_r(r,v)) + r^{2n-3} (2n-2) (-G(r,u)+G(r,v))\notag\\ 
      &=  r^{2n-2} h(r,u,v) (vu'- uv') - r^{2n-3} \int_{v(r)}^{u(r)} rg_r(r,t)+(2n-2) g(r,t)\dt.
    \end{align}
    From $u(r)>v(r)$ for $r>r_0$ and the differential inequality for $g$ from assumption $(A4)$ we obtain that
    the following inequality holds for all $r>r_0$ 
    \begin{equation} \label{Gl Thm2 Ungl}
      \derivord{}{r} \Big( r^{2n-2}\Big(-\frac{1}{2}u'^2+\frac{1}{2}v'^2-G(r,u)+G(r,v) \Big) \Big)
      \geq r^{2n-2} h(r,u,v) (vu'- uv').
	\end{equation}
    From $(A1)$ we get $|g(r,z)| = |f(r,z,0)| \leq C(|z|+|z|^p)$ for all $r,z>0$. In view of~\eqref{Eind Gl
    exp Abf} this entails
    \begin{align*}
      \limsup_{r\to\infty} r^{2n-2}\big| G(r,u(r))-G(r,v(r)) \big|
      &\leq \limsup_{r\to\infty} r^{2n-2} \int_{v(r)}^{u(r)} |g(r,t)|\dt \\
      &\leq C\cdot  \limsup_{r\to\infty}  r^{2n-2} \Big(\frac{1}{2}
      u(r)^2+\frac{1}{p+1} u(r)^{p+1}\Big) \\
      &= 0.         
    \end{align*}
    Hence, integrating the inequality \eqref{Gl Thm2 Ungl} over $(r_0,\infty)$ we obtain
    \begin{align*}
      \int_{r_0}^\infty r^{2n-2}&h(r,u,v)(vu'- uv')\dr \notag\\
      &\leq r_0^{2n-2} \Big( \frac{1}{2} u'(r_0)^2- \frac{1}{2}v'(r_0)^2 + G(r_0,u(r_0))-G(r_0,v(r_0))\Big)
      \notag \\
      &= \frac{1}{2}r_0^{2n-2} (u'(r_0)-v'(r_0))(u'(r_0)+v'(r_0)) \notag \\
      &\leq 0.
    \end{align*}
    In the last inequality we used $(u-v)'(r_0)>0$ and $u'(r_0),v'(r_0)\leq 0$. Notice that the integral on
    the left hand side exists due to $|h(r,u,v)|\leq |f(r,u,v)|+|g(r,u)| \leq 2C(|u|+|v|+|u|^p+|v|^p)$ and
    \eqref{Eind Gl exp Abf}. From
    the inequality $(vu'- uv')(r_0)= u(r_0)(u'-v')(r_0)> 0$ and the positivity of $h$ we infer that there is a point $r_1\in
    (r_0,\infty)$ with $(vu'-uv')(r_1)=0$. We thus obtain
    \begin{align*}
      0
      &= \int_{r_1}^\infty \Big( r^{n-1}(u v' - vu') \Big)'\dr \\
      &=  \int_{r_1}^\infty r^{n-1}\Big( u \Big(v''+\frac{n-1}{r}v'\Big) - v\Big(u''+\frac{n-1}{r}u'\Big)\Big)\dr
      \\
      &= \int_{r_1}^\infty r^{n-1} \big( -uf(r,v,u) + vf(r,u,v) \big) \dr.
    \end{align*}
    Using $(A3')$ we get another point $r_2\in (r_1,\infty)\subset (r_0,\infty)$ with
    $u(r_2)=v(r_2)$ which contradicts the assumption $u-v>0$ on $(r_0,\infty)$. Hence, we get $u=v$. \qed

    \newpage 
    
    \begin{bem} \label{Bem Satz 2}~
      \begin{itemize} 
        \item[(i)] If we drop assumption $(A3)$ then it is not clear to the author whether the function
        $vu'-uv'=v^2(\frac{u}{v})'$ has a zero $r_1$ lying in $(r_0,\infty)$. Without imposing $(A3)$ we
        therefore get the following alternative: Either we have $u=v$ or there is a point $r_0>0$ such that $u(r_0)=v(r_0)$ and
        $\frac{u}{v}$ is strictly monotone on $(r_0,\infty)$.
        \item[(ii)] In the proof we only needed the inequality $(u+v)'(r_0)\leq 0$. As a
        consequence, the assumption that $(u,v)$ is decreasing may be slightly relaxed.
        \item[(iii)] In contrast to Theorem \ref{Eind Thm 1} it is not clear if this result admits
        a generalization that deals with systems of the form \eqref{Gl System QuiSou}.
      \end{itemize}      
    \end{bem}
   
   \section{Proof of Corollary \ref{Eind Cor}}
   
    In order to apply Theorem  \ref{Eind Thm 1} and Theorem \ref{Eind Thm 2} to the system \eqref{Eind Gl DGL}
    we set 
    $$
      f(r,z_1,z_2) = - z_1 + z_1^{2q-1}+bz_1^{q-1}z_2^q
      \qquad\text{for }z_1,z_2,r\geq 0
    $$ 
    where $q$ satisfies \eqref{Gl superlinear subcritical}. The function $f$ clearly satisfies $(A1)$.
    Morover $(A2)$ holds for $m\in (0,1)$ and $\eps>0$ sufficiently small because of the following identities
    for all $r,z_1,z_2>0$ :
    \begin{align*}
      f(r,z_1,z_2)-f(r,z_2,z_1)
      &= -z_1+ z_1^{2q-1}+bz_1^{q-1}z_2^q+z_2-z_2^{2q-1}-bz_2^{q-1}z_1^q \\
      &= (z_1-z_2) \cdot \Big( -1 + \frac{z_1^{2q-1}-z_2^{2q-1}}{z_1-z_2} - b z_1^{q-1}z_2^{q-1} \Big), \\
      f(r,z_1,z_2)+f(r,z_2,z_1)
      &= -z_1 + z_1^{2q-1}+bz_1^{q-1}z_2^q - z_2 + z_2^{2q-1}+bz_2^{q-1}z_1^q) \\
      &= (z_1+z_2) \cdot \Big( -1 + \frac{z_1^{2q-1}+z_2^{2q-1}}{z_1+z_2} + bz_1^{q-1}z_2^{q-1} \Big).
	\end{align*}
	Now let us check for which values of $b,q$ the assumption $(A3)$ or
	$(A3')$ holds true. If ${z_1>z_2>0,r>0}$ and $k:= \frac{z_2}{z_1}\in (0,1)$ then
    \begin{align*}
      f(r,z_1,z_2)z_2-f(r,z_2,z_1)z_1
      &= z_1z_2(z_1^{2q-2}+bz_1^{q-2}z_2^q-z_2^{2q-2}-bz_2^{q-2}z_1^q) \\
      &= z_1^{2q}k (1+bk^q-k^{2q-2}-bk^{q-2})
    \end{align*}
    and we are lead to the study of the function $\psi(k):=1+bk^q-k^{2q-2}-bk^{q-2}$ for
    $k\in (0,1)$. When $1<q<2,b\geq q-1$ or $q=2,b>1$ the function $\psi$ turns out to be negative on the
    interval $(0,1)$ so that Theorem \ref{Eind Thm 1} applies. We obtain $u=v$ and assertion (i) of the
    Corollary follows from Kwong's uniqueness result \cite{Kwo_Uniqueness_of_positive}.
    
    \medskip
    
    Now we apply Theorem \ref{Eind Thm 2} to prove assertion (ii). In case $q=2,b<1$ or $q>2,b\leq q-1$
    the function $\psi$ is positive on $(0,1)$ so that $(A3')$ holds. If moreover $b>0$
    and $n=1$ then $f$ also satisfies $(A4)$. Since $b>0,q\geq 2$ implies that all positive solutions of
    \eqref{Eind Gl DGL} are decreasing (cf. \cite{Iko_Uniqueness_of_positive}, Theorem 4.1) assertion (ii) of
    the Corollary follows from Theorem \ref{Eind Thm 2}. \qed
  
  \section{Appendix}
  
  We show that the uniqueness result found by Ma and Zhao \cite{MaZh_Uniqueness_of_ground} can be derived from
  our Theorem~\ref{Eind Thm 1}. In~\cite{MaZh_Uniqueness_of_ground}, Theorem 2 the authors prove that the positive solution $(u,v)$ of
  the nonlinear Schr\"odinger system
  \begin{align} \label{Gl DGL MaZhao}
     	  \begin{gathered}
      		~\quad  - u'' -\frac{n-1}{r}u' + u =  \mu_1 u^{2q-1} + \beta_1 u^{q-1}v^q
      		\qquad\text{in }\R_{>0},  		\\
      		~\quad - v'' -\frac{n-1}{r}v' + v = \mu_2 v^{2q-1} + \beta_2 u^q v^{q-1}
      		\qquad\text{in }\R_{>0}, \\
      		~\quad u'(0)=v'(0)=0,\quad u,v>0, \quad  u(r),v(r)\to 0 \quad (r\to\infty)
    	  \end{gathered}
  \end{align}
  is unique whenever $q$ satisfies \eqref{Gl superlinear subcritical} and 
  $$
    \mu_1,\mu_2\leq 0,\quad
    \beta_1,\beta_2>0,\quad
    \mu_1\beta_1^{q-1}=\mu_2\beta_2^{q-1},\quad
    \beta_1^{\frac{q-2}{2}}\mu_1 +  \beta_2^{\frac{q}{2}}>0,\quad
    \beta_2^{\frac{q-2}{2}} \mu_2 +\beta_1^{\frac{q}{2}} > 0. 
  $$  
  Now let us show how this result follows from Theorem \ref{Eind Thm 1}. As in
  \cite{MaZh_Uniqueness_of_ground} it suffices to show
  \begin{equation} \label{Eind GL Res MaZhao}
    u = \big( \mu_2 + \beta_1^{\frac{q}{2}}\beta_2^{-\frac{q-2}{2}} \big)^{\frac{1}{2q-2}}
	    \big( \mu_1 +
	    \beta_2^{\frac{q}{2}}\beta_1^{-\frac{q-2}{2}} \big)^{-\frac{1}{2q-2}} v
  \end{equation}  
  for any positive solution $(u,v)$ of \eqref{Gl DGL MaZhao}. A calculation shows that 
  $(\tilde u,\tilde v) = (|\mu_1|^{\frac{1}{2q-2}} u, |\mu_2|^{\frac{1}{2q-2}}v)$ is a solution of
  \begin{align*}
    \begin{gathered}
      		~\quad  - \tilde u'' -\frac{n-1}{r} \tilde u' + \tilde u
      		=  - \tilde u^{2q-1} + \beta \tilde u^{q-1} \tilde v^q \qquad\text{in }\R_{>0},  		\\
      		~\quad - \tilde v'' -\frac{n-1}{r}\tilde v' + \tilde v
      		= - \tilde v^{2q-1} + \beta \tilde u^q \tilde v^{q-1} \qquad\text{in }\R_{>0}, \\
      		~\quad \tilde u'(0)=\tilde v'(0)=0,\quad \tilde u,\tilde v>0, \quad  \tilde u(r), \tilde v(r)\to 0
      		\quad (r\to\infty)
    \end{gathered}
  \end{align*}
  where the parameter $\beta$ is given by 
  $$
    \beta = \beta_1|\mu_1|^{-\frac{q-2}{2q-2}} |\mu_2|^{-\frac{q}{2q-2}} = \beta_2 |\mu_1|^{-\frac{q}{2q-2}}
    |\mu_2|^{-\frac{q-2}{2q-2}}>0.
  $$
  Since the nonlinearity $f(r,z_1,z_2)=-z_1-z_1^{2q-1}+\beta z_1^{q-1}z_2^q$ satisfies the
  assumptions $(A1),(A2)$, $(A3)$ we obtain from Theorem \ref{Eind Thm 1} that $\tilde u$ equals $\tilde v$.
  Due to $\mu_1\beta_1^{q-1}=\mu_2\beta_2^{q-1}$ this implies \eqref{Eind GL Res MaZhao} and that's all we
  had to show.
%

\begin{thebibliography}{11}


\bibitem{BaDaWa_A_Liouville_Theorem}
  Bartsch, T. and Dancer, N. and Wang, Z.-Q, 
  \emph{A {L}iouville theorem, a-priori bounds, and bifurcating branches of positive solutions for a nonlinear
  elliptic system},
  Calc. Var. Partial Differential Equations \textbf{37} (2010), no.~3-4, 345--361.

\bibitem{BusSir_Symmetry_results_for}
  Busca, J. and Sirakov, B., 
  \emph{Symmetry results for semilinear elliptic systems in the whole  space}, 
  J. Differential Equations \textbf{1} (2000), no.~1, 41--56.

\bibitem{Iko_Uniqueness_of_positive}
 Ikoma, N., 
  \emph{Uniqueness of positive solutions for a nonlinear elliptic system}, 
  NoDEA Nonlinear Differential Equations Appl. \textbf{16} (2009), no.~5, 555--567.


\bibitem{Jan_Uniqueness_of_positive}
  Jang, J.,
  \emph{Uniqueness of positive radial solutions of {$\Delta u + f(u)=0$} in {$\mathbb{R}^N$}, {$N\geq 2$}}, 
  Nonlinear Analysis: Theory, Methods and Applications \textbf{73} (2010), no.~7, 2189--2198.


\bibitem{Kwo_Uniqueness_of_positive}
  Kwong, M. K., 
  \emph{Uniqueness of positive solutions of {$\Delta u-u+u^p=0$} in {$\mathbb{R}^n$}}, 
  Arch. Rational Mech. Anal. \textbf{105} (1989), no.~3, 243--266.

\bibitem{MaMoPe_Positive_solutions_for}
  Maia, L. A. and Montefusco, E. and Pellacci, B., 
  \emph{Positive solutions for a weakly coupled nonlinear {S}chr\"odinger system}, 
  J. Differential Equations \textbf{229} (2006), no.~2, 743--767.

\bibitem{MaZh_Uniqueness_of_ground}
  Ma, L. and Zhao, L., 
  \emph{Uniqueness of ground states of some coupled nonlinear
              {S}chr\"odinger systems and their application}, 
  J. Differential Equations \textbf{245} (2008), no.~9, 2551--2565.

\bibitem{PelSer_Uniqueness_of_positive}
   Peletier, L. A. and Serrin, J.,
  \emph{Uniqueness of positive solutions of semilinear equations in
              {$\mathbb{R}^{n}$}}, 
  Arch. Rational Mech. Anal. \textbf{81} (1983), no.~2, 181--197.

\bibitem{QuiSou_Symmetry_of_components}
  Quittner, P. and Souplet, P., 
  \emph{Symmetry of {C}omponents for {S}emilinear {E}lliptic
              {S}ystems}, 
  SIAM J. Math. Anal. \textbf{44} (2012), no.~4, 2545--2559.

\bibitem{Sir_Least_energy_solutions}
  Sirakov, B., 
  \emph{Least energy solitary waves for a system of nonlinear
              {S}chr\"odinger equations in {$\mathbb{R}^n$}}, 
  Comm. Math. Phys. \textbf{271} (2007), no.~1, 199--221.

\bibitem{WeiYao_Uniqueness_of_positive}
  Wei, J. and Yao, W., 
  \emph{Uniqueness of positive solutions to some coupled nonlinear
              {S}chr\"odinger equations}, 
  Commun. Pure Appl. Anal. \textbf{11} (2012), no.~3, 1003--1011.

\end{thebibliography}

\end{document}